\documentclass[11pt]{article}
\usepackage[latin1]{inputenc}
\usepackage[T1]{fontenc}
\usepackage{amsmath}
\usepackage{bbm}
\usepackage{amsthm}
\usepackage{amssymb,mathrsfs}
\usepackage{fourier}
\usepackage[english,french]{babel} 
\usepackage[all]{xy}
\usepackage{setspace}
\usepackage{color}
\usepackage{tabularx}

\newtheorem{theorem}{Theorem}[section]
\newtheorem{corollary}[theorem]{Corollary}

\newtheorem{lemma}[theorem]{Lemma}
\newtheorem{proposition}[theorem]{Proposition}

\newtheorem{definition}[theorem]{Definition}

\allowdisplaybreaks[1]

\newcommand\N{{\mathbb N}}
\newcommand\Z{{\mathbb Z}}
\newcommand\Q{{\mathbb Q}}
\newcommand\R{{\mathbb R}}
\newcommand\C{{\mathbb C}}
\newcommand\T{{\mathbb T}}

\begin{document}

\title{A zero-$\sqrt 5\over 2$ law for cosine families}

\date{\relax}
\author{Jean Esterle}

\maketitle

{\bf Abstract}: Let $a \in \R,$ and let $k(a)$ be the largest constant such that $sup\vert cos(na)-cos(nb)\vert <k(a)$ for $b\in \R$ implies that $b \in \pm a+2\pi\Z. $ We show that
if a cosine sequence $(C(n))_{n\in \Z}$ with values in a Banach algebra $A$ satisfies $sup_{n\ge 1}\Vert C(n) -cos(na).1_A\Vert <k(a),$ then $C(n)=cos(na)$ for $n\in \Z.$ Since
${\sqrt 5\over 2} \le k(a) \le {8\over 3\sqrt 3}$ for every $a \in \R,$ this shows that if some cosine family $(C(g))_{g\in G}$ over an abelian group $G$ in a Banach algebra satisfies $sup_{g\in G}\Vert C(g)-c(g)\Vert <{\sqrt 5\over 2}$ for some scalar cosine family $(c(g))_{g\in G},$ then $C(g)=c(g)$ for $g\in G,$ and the constant ${\sqrt 5\over 2}$ is optimal. 
We also describe the set of all real numbers  $a \in [0,\pi]$ satisfying $k(a)\le {3\over 2}.$

\smallskip


\smallskip


\smallskip

\section{Introduction}

Let $G$ be an abelian group. Recall that a $G$-cosine family of elements of a unital normed algebra $A$ with unit element $1_A$ is a family $(C(g))_{g\in G}$ of elements of $A$  satisfying the so-called d'Alembert equation

\begin{equation}C_0=1_A, C(g+h)+C(g-h)=2C(g)C(h)  \ \ (g\in G,  h\in G).\end{equation}

A $\R$-cosine family is called a cosine function, and a $\Z$-cosine family is called a cosine sequence.

A cosine family $C= (C(g))_{g\in G}$ is said to be bounded if there exists $M>0$ such that $\Vert C(g)\Vert \le M$ for every $g \in G.$ In this case we set

$$\Vert C \Vert_{\infty} =sup_{g \in G}\Vert C(g)\Vert, dist(C_1,C_2)=\Vert C_1-C_2\Vert_{\infty}.$$
A  cosine family is said to be scalar if $C(g)\in \C.1_A$ for every $g \in G.$  It is easy to see and well-known that a bounded scalar cosine sequence satisfies $C(n)=cos(an)$ for some $a \in \R.$  

Strongly continuous operator valued cosine functions are a classical tool in the study of differential equations, see for example \cite{abhn}, \cite{beh}, \cite{n}, \cite{tw}, and a functional calculus approach to these objects was  developped recently in \cite{h1}. 

Bobrowski and Chojnacki proved in \cite{bc} that if a strongly continuous operator valued cosine function on a Banach space $(C(t))_{t\in \R}$ satisfies $sup_{t\ge 0}\Vert C(t) -c(t)\Vert <1/2$ for some scalar bounded continuous cosine function $c(t)$ then $C(t)=c(t)$ pour $t \in \R,$ and Zwart and F. Schwenninger showed in \cite{zs} that this result remains valid under the condition $sup_{t\ge 0}\Vert C(t) -c(t)\Vert <1.$ The proofs were based on rather involved arguments from operator theory and semigroup theory. Very recently, Bobrowski, Chojnacki and Gregosiewicz \cite{bcg} showed more precisely that if a cosine function $C=C(t)$ satisfies sup$_{t\in \R}\Vert C(t)-c(t)\Vert <{8\over 3\sqrt 3}$ for some scalar bounded continuous cosine function $c(t),$ then $C(t)=c(t)$ for $t\in \R,$ without any continuity assumption on $C,$ and the same result was obtained independently by the author in \cite{e0}. The constant ${8\over 3\sqrt 3}$ is obviously optimal, since sup$_{t\in \R}\vert cos(at)-cos(3at)\vert ={8\over 3\sqrt 3}$ for every $a \in \R\setminus \{0\}.$

The author also proved in \cite{e0} that if a cosine sequence $(C(t))_{t\in \R}$ satisfies sup$_{t\in \R}\Vert C(t)-cos(at)1_A\Vert =m<2$ for some $a\neq 0,$ then the closed algebra generated by $(C(t))_{t\in \R}$ is isomorphic to $\C^k$ for some $k\ge 1,$ and that there exists a finite family $p_1,\dots,p_k$ of pairwise orthogonal idempotents of $A$ and a family $(b_1,\dots,b_k)$ of distinct elements of the finite set $\Delta(a,m):=\{ b\ge 0 \ : \ sup_{t\in \R}\vert cos(bt)-cos(at)\vert \le m\}$ such that we have

$$C(t)=\sum \limits _{j=1}^kcos(b_jt)p_j \ \ (j\in \R).$$

Also Chojnacki developped in \cite{c1} an   elementary argument to show that if $(C(n))_{n\in \Z}$ is a cosine sequence in a unital normed algebra $A$ satisfying $sup_{n\ge 1}\Vert C(n)-c(n)\Vert <1$ for some scalar cosine sequence $(c(n))_{n\in Z}$ then $c(n)=C(n)$ for every $n,$
which obviously implies the result of Zwart and F. Schwenninger. His approach is based on an elaborated adaptation of a very short  elementary argument used by Wallen in \cite {w} to prove an improvement of the classical Cox-Nakamura-Yoshida-Hirschfeld-Wallen
theorem \cite{co}, \cite{h}, \cite{ny} which shows that if an element $a$ of a unital normed algebra $A$ satisfies $sup_{n\ge 1}\Vert a^n-1\Vert <1,$ then $a=1.$

Applying this result to the cosine sequences $C(ng)$ and $c(ng)$ for $g\in G,$  Chonajcki observed in \cite{c1} that if a cosine family $C(g)$  satisfies $sup_{g\in G}\Vert C(g)-c(g)\Vert
 <1$ for some scalar cosine family $c(g)$ then $C(g)=c(g)$ for every $g\in G.$
 
 In the same direction Schwenninger and Zwart showed in \cite{zs0} that if a cosine sequence $(C(n))_{n\in \Z}$ in a Banach algebra $A$ satisfies $sup_{n\ge 1}\Vert C(n)-1_A\Vert <{3\over 2},$ then $C(n)=1_A$ for every $n.$
 
 The purpose of this paper is to obtain optimal results of this type. We prove a "zero-${\sqrt 5\over 2}$" law :  if a cosine family $(C(g))_{g\in G}$  satisfies $sup_{g\in G}\Vert C(g)-c(g)\Vert
 <{\sqrt 5\over 2}$ for some scalar cosine family $(c(g))_{g \in G}$ then $C(g)=c(g)$ for every $g \in G.$ Since sup$_{n\ge 1}\left \vert cos \left ({2n\pi\over 5}\right )-cos\left ({4n\pi\over 5}\right ) \right \vert=cos\left ({2\pi\over 5}\right )-cos\left ({\pi\over 5}\right )={\sqrt 5\over 2},$ the constant ${\sqrt 5\over 2}$ is optimal. 
 
 In fact for every $a \in \R$ there exists a largest constant $k(a)$
 such that sup$_{n\ge 1}\vert cos(nb)-cos(na)\vert <k(a)$ implies that $cos(nb)=cos(na)$ for $n\ge 1,$ and we prove that if a cosine sequence $(C(n))_{n\in \Z}$ in a Banach algebra $A$ satisfies sup$_{n\ge 1}\vert C(n)-cos(na)1_A\vert <k(a)$ then $C(n)=cos(na)$ for $n\ge 1.$ This follows from the following result, proved by the author in \cite{e0}.
 
 \begin{theorem} Let $(C(n))_{n\in \Z}$ be a bounded cosine sequence in a Banach algebra. If $spec(C(1))$ is a singleton, then the sequence $(C(n))_{n\in \Z}$ is scalar, and so there exists $a\in \R$ such that $C(n)=cos(na)$ for $n\ge 1.$
 \end{theorem}

 The second part of the paper is devoted to a discussion of the values of the constant $k(a).$ As mentioned above, it follows from \cite{zs0} that $k(0)={3\over 2},$ and it is obvious that $k(a)\le $sup$_{n\ge 1}\vert cos(na)-cos(3na)\vert \le {8\over 3\sqrt 3}$ if $a \notin {\pi\over 2}\Z.$ We observe that $k(a)={8\over 3\sqrt 3}$ if ${a\over \pi}$ is irrational, and we prove, using basic results about cyclotomic fields, that  $k(a)<{8\over 3\sqrt 3}$ if ${a\over \pi}$ is rational.
 
 We also show that the set $\Omega(m):=\{ a \in [0, \pi] \ \vert \ k(a)\le m\}$ is finite for every $m<{8\over 3\sqrt 3}.$ We describe in detail the set $\Omega\left ({3\over 2}\right):$ it contains 43 elements, and the only values for $k(a)$ for which 
 $k(a) <{3\over 2}$ are ${\sqrt 2 \over 5}=cos\left ( {\pi\over 5}\right )+cos\left ( {2\pi\over 5}\right )\approx 1.1180,$ $\sqrt 2=cos\left ( {\pi\over 4}\right )+cos\left ( {3\pi\over 4}\right )\approx 1.4142,$ and $cos\left ( {2\pi\over 11}\right )+cos\left ( {3\pi\over 11}\right )\approx 1.4961.$
 
  The zero-${\sqrt 5\over 2}$ law follows then from the fact that $k(a)\ge cos\left ({\pi\over 5}\right ) + cos\left ({\pi\over 5}\right )={\sqrt 5\over 2}$ for every $a\in \R.$

We also show that given $a \in \R$ and $m <2$ the set $\Gamma(a,m)$ of scalar cosine sequences $(c(n))_{n\in \Z}$ satisfying $sup_{n\in \Z}\vert c(n) -cos(na)\vert \le m$ is finite. This implies that if a cosine sequence $(C(n))_{n\in \Z}$ satisfies sup$_{n\in \Z}\Vert C(n)-cos(an)1_A\Vert \le m,$ then there exists $k\le card\left (\Gamma(a,m)\right )$ such that the closed algebra generated by $(C(n))_{n\in \Z}$ is isomorphic to $\C^k,$ and there exists a finite family $p_1, \dots, p_k$ of pairwise orthogonal idempotents of $A$ and a finite family 
$c_1, \dots, c_k$ of distinct elements of $\Gamma(a,m)$ such that we have

$$C(n)=\sum \limits _{j=1}^kc_j(n)p_j \ \ (n \in \Z).$$

This last result does not extend to cosine families over general abelian group. Let $G=(\Z/3\Z)^\N$: we give an easy example of a $G$-cosine family $(C(g))_{g\in G}$ with values in $l^{\infty}$ such that the closed subalgebra generated by $(C(g))_{g\in G}$ equals $l^{\infty},$ while  sup$_{g\in G}\Vert 1_{l^{\infty}} -C(g)\Vert ={3\over 2}.$

The author warmly thanks Christine Bachoc and Pierre Parent for giving him the arguments from number theory which lead to a simple proof of the fact that $k(a)<{8\over 3\sqrt 3}$ if $a\notin \pi\Q.$
 \section{Distance between bounded scalar cosine sequences}

We introduce the following notation, to be used throughout the paper.

\begin{definition} Let $a \in \pi \Q.$ The order of $a$, denoted by $ord(a),$ is the smallest integer $u\ge 1$ such that $e^{iua}=1.$
\end{definition}

Recall that a subset $S$ of the unit circle $\T$ is said to be independent if $z_1^{n_1}\dots z_k^{n^k}\neq 1$ for every finite family $(z_1,\dots,z_k)$ of distincts elements of $S$ and every family $(n_1, \dots,n_k)\in \Z^k$ of such that $z_j\neq 0$ for $j\le j\le k.$ It follows from a classical theorem of Kronecker, see for example \cite{ks}, page 21 that if $S=\{z_1,\dots,z_k\}$ is a finite independent set then the sequence $( z_1^n,\dots,z_k^n)_{ n\ge 1}$ is dense in $\T^k.$ We deduce from Kronecker's theorem the following observation.

\begin{proposition} Let $a \in [0,\pi].$ For $m\ge 0,$ set 
$$\Gamma(a,m)=\left \{ b \in [0,\pi] :sup_{n\ge 1} \left | cos(na)-cos(nb)\right | \le m \right \}.$$  Then $\Gamma(a,m)$ is finite for every $m<2.$
\end{proposition}

Proof: Fix $m \in [1,2).$ Notice that if $b\in \R,$ and if the set $\left \{e^{ia},e^{ib}\right \}$ is independent, then it follows from Kronecker's theorem that the sequence $\left (\left (e^{ina},e^{inb}\right )\right )_{n\ge 1}$ is dense in $\T^2,$
and so $sup_{n\ge 1} \left | cos(na)-cos(nb)\right |=2,$ and $b\notin \Gamma(a,m).$

Suppose that ${a\over \pi}\in \Q,$ and denote by $u$ the order of $a,$ so that $e^{iua}=1.$  If ${b\over \pi}\notin \Q,$ then the sequence $\left (e^{iunb} \right )_{n\ge 1}$ is dense in $\T,$ and so
$$2\ge sup_{n\ge 1}\left | cos(na)-cos(nb)\right |\ge  sup_{n\ge 1} \left | 1-cos(nub)\right |=2,$$

which shows that $b\notin \Gamma(a,m).$

The same argument shows that if ${a\over \pi} \notin \Q,$ and if ${b\over \pi} \in \Q,$ then $b\notin \Gamma(a,m).$ So we are left with two situations

1) ${a\over \pi} \notin Q,$ and there exists $p\neq 0,$ $q \neq 0$ and $k\in \Z$ such that $bq=ap+2k\pi.$

\smallskip

2) ${a\over \pi}\in \Q$ and ${b\over \pi} \in \Q.$ 

We consider the first case. Replacing $b \in [0,\pi]$ by $-b\in [-\pi,0]$ if necessary we can assume that $p\ge 1$ and $q\ge 1,$ and we can assume that we have

$$qb =pa +{2k\pi\over r},$$

with $gcd(p,q)=1, r\ge 1, gcd(r,k)=1$ if $k\neq 0.$

We have, since ${ra\over \pi} \notin \Q,$

$$sup_{n\ge 1}\left | cos(na)-cos(nb)\right | \ge sup_{n\ge 1}\left | cos(nrqa)-cos(nrqb)\right | $$ $$= 
sup_{n\ge 1}\left | cos(nrqa)-cos(nrpa)\right | =sup_{t \in \R}\left | cos(qt)-cos(pt) \right |,$$

Since $gcd(p,q)=1,$ we have $sup_{t \in \R}\left | cos(qt)-cos(pt) \right |=2$ if $p$ or $q$ is even, so we can assume that $p$ and $q$ are odd. Set $s={q-1\over 2}.$

It follows from Bezout's identity that there exist $n\ge 1$ such that $e^{2inp\pi\over q}=e^{2is\pi \over q}$ and setting $t={2n\pi\over q},$ we obtain

$$sup_{t \in \R}\left | cos(qt)-cos(pt) \right | \ge 1 -cos\left ({2s\pi\over 2s+1}\right )=1+cos\left ({\pi\over q}\right ).$$

The same argument shows that we have

$$sup_{t \in \R}\left | cos(qt)-cos(pt) \right |\ge 1+cos\left ({\pi\over p}\right ).$$

We obtain

$$p \le {\pi\over arccos(m-1)}, q \le {\pi\over arccos(m-1)}.$$

We also have

$$sup_{n\ge 1}\left | cos(na)-cos(nb)\right | \ge sup_{n\ge 1}\left | cos(nqa)-cos(nqb)\right | $$ $$=sup_{n\ge 1}\left | cos(nqa)-cos\left (npa +{2nkq\pi\over r}\right )\right |.$$

Assume that $k\neq 0,$ set $d=gcd(r,q), r_1={r\over d}, q_1={q\over d}.$ Then $gcd(kq_1,r_1)=1,$ and so there exists $u \ge 1$ such that ${2ukq_1\pi\over r_1}\in {2\pi\over r_1}+2\pi \Z.$ This gives

$$sup_{n\ge 1}\left | cos(na)-cos(nb)\right | \ge sup_{n\ge 1}\left | cos(nuqa)-cos\left (npua +{2n\pi\over r_1}\right )\right |.$$

If $r_1$ is even, set $r_2={r_1\over 2}.$ We obtain

$$ sup_{n\ge 1}\left | cos(nuqa)-cos\left (npua +{2n\pi\over r_1}\right )\right |$$ 

$$\ge sup_{n\ge 0}\left | cos((2n+1)r_2uqa) -cos\left ((2n+1)r_2upa) +\pi \right)\right |.$$

Since $2r_2ua\notin \pi\Q,$ there exists a sequence $(n_j)_{j\ge 1}$ of integers such that $$lim_{j\to +\infty}\left | e^{i2n_jr_2ua+ir_2ua}\right |=1,$$ so that 
$$lim_{j\to +\infty}\left | cos((2n_j+1)r_2uqa) -cos\left ((2n_j+1)r_2upa) +\pi \right)\right |=2,$$ and in this situation $sup_{n\ge 1}\left | cos(na)-cos(nb)\right | =2.$

So we can assume that $r_1$ is odd. Set $r_2={r-1\over 2}.$ The same calculation as above gives
$$ sup_{n\ge 1}\left | cos(nuqa)-cos\left (npua +{2n\pi\over r_1}\right )\right |$$ $$\ge sup_{n\ge 1} \left  | cos((n(2r_2+1)+r_2)uqa) -cos\left ((n(2r_2+1)+r_2)upa +{2(n(2r_2+1)+r_2)\over 2r_2+1}\pi \right)\right | $$ $$\ge1+  cos \left ({\pi\over 2r_2+1}\right ).$$

Hence $r_1=2r_2+1\le {\pi\over arccos(m-1)},$ $r=r_1d\le r_1q\le \left ({\pi\over arccos(m-1)}\right )^2.$

This gives

$$2\vert k \vert \pi \le r \vert qb -pa\vert \le 2\pi \left ({\pi\over arccos(m-1)}\right )^3, \vert k \vert \le  \left ({\pi\over arccos(m-1)}\right )^3.$$

We see that $\Gamma(a,m)$ is finite if ${a\over \pi}\notin Q,$ and that we have

$$card\left ( \Gamma(a,m)\right ) \le  2\left ({\pi\over arccos(m-1)}\right )^7.$$

Now consider the case where ${a\over \pi} \in \Q,$ ${b\over \pi} \in \Q.$ We first discuss the case where $a=0, b\neq 0.$ We have $b={p\pi\over q},$ where $1\le p \le q,$ $gcd(p,q)=1$

If $p=q=1,$ then $b=\pi,$ and $sup_{n\ge 1}\vert 1-cos(n\pi)\vert =2.$ So we may assume that $p\le q-1.$ If $p$ is odd, then we have

$sup_{n\ge 1}\vert 1 -cos(nb)\vert \ge \vert 1-cos(qb)\vert =1-cos(p\pi)=2.$

So we can assume that $p$ is even, so that $q$ is odd. Set $r={q-1\over 2}.$ There exists $n_0\ge 1$ and $r \in \Z$ such that $n_0p -r \in q\Z,$ and we have

$$sup_{n\in \Z}\vert 1 - cos(nb)\vert \ge \left \vert 1  - cos(2n_0b) \right \vert =\left \vert 1 -cos\left ({2r\pi\over 2r+1}\right ) \right | =1+cos\left ({\pi\over q}\right ).$$

We obtain again $q \le {\pi\over arccos(m-1)},$ and $card \left ( \Gamma(0,m)\right )\le \left ( {\pi\over arccos(m-1)}\right )^2.$

Now assume that $a\neq 0,$ and let $u\ge 2$ be the order of $a.$ We have

$$sup_{n\ge 1} \vert 1 -cos(nub)\vert =sup_{n\ge 1} \vert cos(nua) -cos(nub)\vert \le m,$$ and so there exists  there exists $c\in \Gamma(0,m)$ such that $cos(nc)=cos(nub)$ for
$n\ge 1.$ In particular $cos( c )=cos(ub),$ and $b=\pm{c\over u}+{2k\pi\over u},$ where $k \in \Z.$ We obtain

$$card\left (\Gamma(a,m)\right ) \le 2u card\left ( \Gamma(0,m)\right )\le 2ucard\left ({\pi\over arccos(m-1)}\right )^2.$$
$\square$

We do not know whether it is possible to obtain a majorant for $card\left (\Gamma(a,m)\right )$ which depends only on $m$ when $a\in \pi\Q.$ 

\begin{theorem} Let $a\in \R,$ let $m<2,$ and let $(C(n))_{n\in \Z}$ be a cosine sequence in a Banach algebra $A$ such that $sup_{n\ge 1}\Vert C(n)-cos(na)\Vert \le m.$ Then  there exists $k\le card\left (\Gamma(a,m)\right )$ such that the closed algebra generated by $(C(n))_{n\in \Z}$ is isomorphic to $\C^k,$ and there exists a finite family $p_1, \dots, p_k$ of pairwise orthogonal idempotents of $A$ and a finite family  
$b_1, \dots, b_k$ of distinct elements of $\Gamma(a,m)$ such that we have

$$C(n)=\sum \limits _{j=1}^kcos(nb_j)p_j \ \ (n \in \Z).$$

\end{theorem}


Proof: Since $c_n =P_n(c_1),$ where $P_n$ denotes the $n$-th Tchebishev polynomial, $A_1$ is the closed unital subalgebra generated by $c_1$ and the map $\chi \to \chi(c_1)$ is a bijection from $\widehat{A_1}$ onto $spec_{A_1}(c_1).$ Now let $\chi \in \widehat{A_1}.$ The sequence $(\chi(c_n))_{n\ge 1}$ is a scalar cosine sequence, and we have

$$sup_{n\ge 1}\left \vert cos(na)-\chi(c_n)\right \vert <2.$$

It follows then from proposition $2.2$ that $spec_{A_1}(c_1):=\left \{ \lambda = \chi(c_1) : \chi \in \widehat{A_1}\right \}$ is finite. Hence $\widehat{A_1}$ is finite. Let $\chi_1, \dots, \chi_m$ be the elements of $\widehat{A_1}.$ It follows from the standard one-variable holomorphic functional calculus, se for example \cite{da}, that there exists for every $j\le m$ an idempotent $p_j$ of $A_1$ such that $\chi_j(p_j)=1$ and $\chi_k(p_j)=0$ for $k\neq j.$ Hence $p_jp_k=0$ for $j \neq k,$ and $\sum \limits_{j=1}^mp_j$ is the unit element of $A_1.$

Let $x\in A_1.$ Then $(p_jc_n)_{n\in \Z}$ is a cosine sequence in the commutative unital Banach algebra $p_jA_1,$ and $spec_{p_jA_1}(p_jc_1)=\{\chi_j(c_1)\}.$ 

Since $sup_{n\ge 1}\left \Vert p_jcos(na)-p_jc_n \right \Vert \le 2 \Vert p_j\Vert,$ the sequence $(p_jc_n)_{n\ge 1}$ is bounded, and it follows from theorem 2.3 that $(p_jc_n)_{n\ge 1}$ is a scalar sequence, and there exists $\beta_j\in [0,\pi]$ such that $p_jc_n= \chi_j(c_n)p_j= cos(n\beta_j)p_j$ for $n \in \Z.$

Hence $c_n=\sum \limits_{j=1}^m\chi_j(c_n)p_j=\sum \limits_{j=1}^mcos(n\beta_j)p_j$ for $n \ge 1.$ Since $A_1$ is the closed subalgebra of $A$ generated by $c_1,$ we have
$x= \sum \limits_{j=1}^m\chi_j(x)p_j$ for every $x \in A_1,$ which shows that $A_1$ is isomorphic to $\C^m.$ $\square$

\begin{corollary} Let $a\ge 0\in \R,$ and let $k(a)$ be the largest positive real number $m$ such that $\Gamma(a,m)=\{a\}$ for every $m <k(a).$ If $(C(n))_{n\in \Z}$ is a cosine sequence in a Banach algebra $A$ such that $sup_{n\ge 1}\Vert C(n)-cos(na)1_A\Vert <k(a),$ then $C(n)=cos(na)1_A$ for $n\in \Z.$
\end{corollary}

Theorem 2.3 does not extend to cosine  families over general abelian groups, as shown by the following easy result.

\begin{proposition} Let $G:=(\Z/3\Z)^{\N}.$ Then there exists a $G$-cosine family $(C(g))_{g\in G}$ with values in $l^{\infty}$ which satisfies the two following conditions

\smallskip

(i) $sup_{g\in G}\Vert 1_{l^{\infty}}-C(g)\Vert ={3\over 2},$

\smallskip

(ii) The algebra generated by the family $(C(g))_{g\in G}$ is dense in $l^{\infty}.$

\end{proposition}

Proof: Elements $g$ of $G$ can be written under the form $g=(\overline g_m)_{m\ge 1},$ where $g_m\in \{ 0, 1, 2\}.$ Set

$$C(g):= \left (cos\left ( {2g_m\pi\over 3}\right )\right )_{m\ge 1}.$$

Then $(C(g))_{g\in G}$ is a $G$-cosine family with values in $l^{\infty}$ which obviously satisfies (i) since $cos\left ( {2\pi\over 3}\right )= cos\left ( {4\pi\over 3}\right )=-{1\over 2}.$

Now let $\phi=(\phi_m)_{m\in \Z}$ be an idempotent of $l^{\infty},$ and let $S:=\{ m\ge 1 \ \vert \ \phi_m=1\}.$ Set $g_m=1$ if $m\in S,$ $g_m=0$ if $m\ge 1, m\notin S,$ and set $g=(\overline g_m)_{m\ge 1}.$ We have

$$C(0_G)-C(g)=1_{l^{\infty}}-C(g)={3\over 2}\phi,$$

and so $\phi \in A.$ We can identify $l^{\infty}$ to $\mathcal C(\beta\N),$ the algebra of continuous functions on the Stone-C${\breve{\mbox e}}$ch compactification of $\N,$ and $\beta N$ is an extremely disconnected compact set, which means that the closure of every open set is open, see for example \cite{ap}, chap. 6, sec. 6. Since the characteristic function of every
open and closed subset of $\beta \N$ is an idempotent of $l^{\infty},$ the idempotents of $l^{\infty}$ separate points of $\beta \N,$ and it follows from the Stone-Weierstrass theorem
that $A$ is dense in $l^{\infty},$ which proves (ii). $\square$

\section{The values of the constant $k(a)$}

It was shown in \cite{zs0} that $k(0)={3\over 2}.$ We also have the following result.

\begin{proposition} We have $k(a)={8\over 3\sqrt 3}$  if $a\over\pi$ is irrrational, and  $k(a)<{8\over 3\sqrt 3}$  if $a\over\pi$ is rational .
\end{proposition}

Proof: Assume that  ${a\over \pi} \notin \Q.$  Then $3a \notin \pm a +2\pi \Z,$ and we have

$$k(a) \le sup_{n\ge 1}\vert cos(na)-cos(3na)\vert = sup_{x\in \R}\vert cos(x)-cos(3x)\vert={8\over 3\sqrt 3}.$$

We saw above that If ${b\over \pi}$ in $\Q,$ then sup$_{n\ge 1}\vert cos(na)-cos(nb)\vert=2,$ and we also have sup$_{n\ge 1}\vert cos(na)-cos(nb)\vert=2$ if $pa -qb\notin 2\pi\Z$ for $(p,q)\neq (0,0).$ So if sup$_{n\ge 1}\vert cos(na)-cos(nb)\vert <2,$ there exists $p\in Z\setminus\{0\}$, $q \in \Z\setminus\{0\}$ and $r\in \Z$ such that
$pa -qb=2r\pi.$

If $p\neq \pm q$ then it follows from lemma 3.5 of \cite{e0} that we have

$$\mbox{sup}_{n\ge 1}\vert cos(na)-cos(nb)\vert \ge \mbox{sup}_{n\ge 1}\vert cos(nqa)-cos(nqb)\vert =\mbox{sup}_{n\ge 1}\vert cos(qna)-cos(pna)\vert$$
$$=\mbox{sup}_{x\in \R}\vert cos(qx)-cos(px)\vert=\mbox{sup}_{x\in \R}\left \vert cos\left ({p\over q}x\right )-cos(x)\right \vert \ge {8\over 3\sqrt 3}.$$

We are left with the case where $b=\pm a +{2s\pi\over r},$ where $r \in \Z \setminus \{-1,0,1\},$ and we can restrict attention to the case where $b=a+{2s\pi\over r}$ where $r \ge 2,$
$1\le s \le r-1,$ $gcd(r,s)=1.$ It follows from Bezout's identity that there exists for every $p\ge 1$ some $u \in \Z$ such that $ub -ua -{2p\pi\over r} \in 2\pi \Z.$  If $r$ is even, set $p={r\over 2}.$ We have, since the set $\left \{e^{i(2n+1) a}\right \}_{n\ge 1}$ is dense in the unit circle,

$$sup_{n\ge 1}\vert cos(nb)-cos(na)\vert =sup_{n\in \Z}\vert cos(nb)-cos(na)\vert  $$ $$\ge sup_{n\ge 1}\vert cos((2n+1)u b)-cos((2n+1)u a)\vert$$ $$=2sup_{n\ge 1}\vert cos((2n+1)u a)\vert=2.$$

Now assume that $r$ is odd, and set $p={r-1\over 2}.$ We have

$$sup_{n\ge 1}\vert cos(nb)-cos(na)\vert \ge sup_{n\ge 1}\vert cos((2n+1)u b)-cos((2n+1)ua)\vert$$ $$\ge sup_{n\ge 1} \left \vert cos((2nr+1)ua) -cos\left ((2nr+1)u a +(2nr+1)\left (\pi -{\pi\over r}\right )\right )\right \vert $$ $$\ge sup_{x\in \R}\left \vert cos(x) +cos\left(x-{\pi\over r}\right )\right | \ge 2cos \left ({\pi\over 2r}\right )\ge \sqrt 3 >{8\over 3\sqrt 3}.$$ 

Now assume that $a\over \pi$ is rational. If the order of $a$ is equal to $1,$ then $k(a)=1.5,$ and we will see later that this is also true if the order of $a$ equals $2$ or $4.$

Otherwise we have

$$k(a) \le sup_{n\ge 1}\vert cos(na)-cos(3na)\vert=max_{1\le n \le u}\vert cos(na)-cos(3na)\vert.$$

We have $\vert cos(nx) -cos(3nx)\vert <{8\over \pi \sqrt 3}$ if $x \notin \pm arccos\left ({1\over \sqrt 3}\right )+\pi\Z.$ If $na\in \pm arccos\left ({1\over \sqrt 3}\right )+\pi\Z$ for some $n\ge 1,$ then ${ arccos\left ({1\over \sqrt 3}\right) \over \pi}$ would be rational, and $\alpha:= {1\over \sqrt 3} +{\sqrt 2 i\over \sqrt 3}$ would be a root of unity. So $\beta =\alpha^2=-{1\over 3}+{2\sqrt 2i\over 3}$ would have the form $\beta =e^{2ik\pi\over n}$ for some $n\le 1$ and some positive integer $k\ge n$ such that $gcd(k,n)=1.$

Let $\Q(\beta)$ be the smallest subfield of $\C$ containing $\Q \cup {\beta}.$ Since $3\beta^2 +2\beta +3=0,$ the degree of $\Q(\beta)$ over $\Q$ is equal to $2.$ On the other hand  the Galois group $Gal(\Q(\beta)/\Q)$ is isomorphic to $\left (\Z/ n\Z\right )^{\times},$ the group of invertible elements of $\Z/n\Z$, and we have, see \cite{wa}, theorem 2.5

$$H(n)= deg(\Q(\beta)/\Q)=2,$$

where $H(n)=card\left (\left (\Z/ n\Z\right )^{\times}\right )$ denotes the number of integers $p\in \{1,\dots,n\}$ such that $gcd(p,n)=1.$

Let $P(n)$ be the set of prime divisors of $n.$ It is weil-known that we have, writing $n=\Pi_{p\in P(n)}p^{\alpha_p},$ see for example \cite{wa}, exercise 1.1,

$$H(n)=\Pi_{p\in P(n)}p^{\alpha_p-1}(p-1).$$

It follows immediately from this identity that the only possibilities to get $H(n)=2$ are $n=3,$ $n=4,$ and $n=6.$ Since $\beta^3\neq 1,$ $\beta^4\neq 1,$ and $\beta^6\neq 1,$ we see
that $\beta\over \pi$ is irrational, and so $k(a)<{8\over 3\sqrt 3}$ if ${a\over \pi}$ is rational. $\square$

$\square$

We know that if  ${a\over \pi}$ is rational, and if ${b\over \pi}$ is irrational, then $sup_{n\ge 1}\vert cos(na)-cos(nb)\vert =2.$  We discuss now the case where ${a\over \pi}$ and
 ${b\over \pi}$ are both rational, with $b\notin \pm a+ 2\pi \Z.$

\begin{lemma} Let $a,b\in (0,\pi].$ 

(i) If  $7a\le b \le {\pi\over 2},$ or if ${\pi\over 2} \le b \le {5\pi\over 6},$ with $\left \vert b-{2\pi\over 3}\right \vert \ge 7a,$ then

$$sup_{n\ge 1}\left \vert cos \left ({na}\right )-cos \left ({nb}\right )\right \vert > 1.55.$$

(ii) If ${5\pi\over 6}\le b \le \pi ,$ and if $b \ge 4a,$ then 

$$cos(a)-cos(b)>1.57.$$

\end{lemma}

Proof: (i) Assume that $7a\le b \le {\pi\over 2},$ let $p$ be the largest integer such that $pb< {3\pi \over 4},$ and set $q=p+1.$ We have ${3\pi\over 4} \le qb\le {5\pi\over 4},$ $0\le qa \le {5\pi\over 28},$ and we obtain

$$sup_{n\ge 1} \left \vert cos \left ({na}\right )-cos \left ({nb}\right )\right \vert\ge  cos \left ({qa}\right )-cos \left ({qb}\right )
\ge cos\left ({5\pi \over 28}\right)  +cos\left ({\pi\over 4}\right )>1.55.$$

Now assume that  ${\pi\over 2} \le b \le {5\pi\over 6},$ with $\vert b-{2\pi\over 3}\vert \ge 7a,$ and set $c=\vert 3b-2\pi \vert.$ Since $\left \vert b -{2\pi\over 3}\right \vert \le {\pi\over 6},$ we have $21a\le c \le {\pi\over 2},$ and we obtain

$$sup_{n\ge 1}\left \vert cos(na) -cos(nb)\right \vert \ge sup_{n\ge 1}\left \vert cos(3na) -cos(3nb)\right \vert $$ $$= sup_{n\ge 1}\left \vert cos(3na) -cos(nc)\right \vert >1.55.$$

(ii) If ${5\pi\over 6}\le b \le \pi ,$ and if $b \ge 4a,$ then $0<a \le {\pi\over 4},$ and we have

$$cos(a)-cos(b) \ge cos\left ( {\pi\over 4}\right ) +cos\left ( {\pi\over 6}\right )>1.57.$$

\begin{lemma}  Let $p, q$ be two positive integers such that $p<q.$ 

(i) If $q \neq 3p,$  then there exists $u_{p,q}\ge 1$ such that, if $ord(a)\ge u_{p,q}$ we have

$$sup_{n\ge 1}\vert cos(npa)-cos(nqa)\vert >{8\over \sqrt 3}.$$

(ii) If $q=3p,$ then for every $m<{8\over 3\sqrt 3}$ there exists $u_p(m) \ge 1$ such that if $ord(a)\ge u(m)$ we have

$$sup_{n\ge 1}\vert cos(npa)-cos(3npa)\vert >m.$$

\end{lemma}

Proof: Set $\lambda=sup_{x\in \R}\vert cos(px)-cos(qx)\vert =sup_{x\ge 0}\vert cos(px)-cos(qx)\vert.$ An elementary verification shows that $\lambda>{8\over 3\sqrt 3}$ if $q\neq 3p,$
and $\lambda ={8\over 3\sqrt 3}$ if $q=3p,$ see for example \cite{e0}. Now let $\mu <\lambda,$ and let $\eta<\delta$ be two real numbers such that $\vert cos(px)-cos(qx)\vert >\mu$ for 
$\eta\le x \le \delta.$ Since $\{e^{ian}\}_{n\ge 1}=\{e^{2ni\pi\over u}\}_{1\le n\le u},$ we see that $sup_{n\ge 1}\vert cos(npa)-cos(nqa)\vert >\mu$ if ${2\pi\over u}<\delta - \eta,$ and the lemma follows. $\square$

\begin{lemma} Assume that ${a\over \pi}$ and ${b\over \pi}$ are rational, let $u\ge 1$ be the order of $a$ and let $v$ be the order of $b.$ 

(i) If $u \neq v,$ $u\neq 3v,$ $v \neq 3u$ then $sup_{n\ge 1}\vert cos(na) -cos(nb)\vert \ge 1+cos\left ( {\pi\over 5}\right )>1.8 > {8\over 3\sqrt 3}.$

(ii) If $u=v,$ and if $b \notin \pm a +2\pi\Z,$ then there exists $w \in \Z$ such that $2\le w\le {u\over 2}$ and $gcd(u,w)=1$ satisfying \begin{equation}sup_{n\ge 1}\cos(na)-cos(nb)\vert = sup_{n\ge 1} \left \vert cos\left ({2 n\pi\over u}\right) -cos \left ( {2nw\pi\over u}\right ) \right \vert.\end{equation}

Conversely if $a\in \pi\Q$ has order $u,$ then for every integer $w$ such that $gcd(w,u)=1,$ there exists $b\in \pi \Q$ of order $u$ satisfying (2).

(iii) If $v=3u,$ then there exists an integer $w$ such that $1\le w\le {u\over 2}$ and $gcd(u,w)=1$ satisfying \begin{equation}sup_{n\ge 1}\cos(na)-cos(nb)\vert = sup_{n\ge 1} \left \vert cos\left ({2 n\pi\over 3u}\right) -cos \left ( {2nw\pi\over u}\right ) \right \vert.\end{equation} 

Conversely if $a\in \pi\Q$ has order $u,$  then for every integer $w$ such that $gcd(w,u)=1$ there exists $b\in \pi\Q$ of order $3u$ satisfying (3).
   
   (iv) If $u=3v,$ then there exists an integer $w$ such that $1\le w\le {u\over 6}$ and $gcd\left ({u\over 3},w\right )=1$ satisfying \begin{equation}sup_{n\ge 1}\cos(na)-cos(nb)\vert = sup_{n\ge 1} \left \vert cos\left ({2 n\pi\over u}\right) -cos \left ( {6nw\pi\over u}\right ) \right \vert.\end{equation} 
   
Conversely  if the order $u$ of $a\in \pi\Q$ is divisible by $3,$  then for every integer $w$ such that $gcd\left ({u\over 3},w\right )=1$ there exists $b\in \pi \Q$ of order $u\over 3$ satisfying (4).

\end{lemma}

Proof: (i) Assume that $u \neq v,$ say, $u<v,$ and let $w \neq 1$ be the order of $ub,$ which is a divisor of $v.$ We have $ub={2\pi \alpha\over w},$ with $gcd(\alpha,w)=1,$ and there exists $\gamma \ge 1$ such that $\alpha \gamma -1 \in w\Z.$ We obtain

$$sup_{n\ge 1}\vert cos(na) -cos(nb) \vert \ge sup_{n\ge 1} \vert cos(nu\gamma a) -cos(nu\gamma b)\vert =sup_{1\le n \le w } 1-cos\left ({2n\pi\over w}\right ).$$ 

 If $w$ is even, then  $sup_{n\ge 1}\vert cos(na) -cos(nb) \vert=2.$ If $w$ is odd, set $s={w-1\over 2}.$ We obtain
 
 $$sup_{n\ge 1}\vert cos(na) -cos(nb) \vert \ge 1 -cos\left ({2s\pi\over w}\right ) = 1 +cos\left ( {\pi\over w}\right ).$$
 
 If $w\ge 5,$ we obtain
 
 $$sup_{n\ge 1}\vert cos(na) -cos(nb) \vert \ge 1+cos\left ({\pi\over 5}\right )>1.8>{8\over 3\sqrt 3}.$$
 
  If $w=3,$ let $d=gcd(u,v),$ and set $r={u\over d}.$ Then $w=3={v\over d}>r.$  So either $r=1$ or $r=2.$
  
  If $r=2,$ we have $u=2d, v=3d,$ $a={2p\pi\over 2d}={p\pi\over d}$ with $p$ odd, $b={2q\pi\over 3d}$ with $gcd(q,3d)=1,$ and we obtain
  
   $$sup_{n\ge 1}\vert cos(na) -cos(nb) \vert = \left \vert cos\left ( {3da}\right ) -cos\left ( {3db}\right ) \right \vert $$ $$\ge \vert cos(3p\pi)-cos(2q\pi)\vert =2.$$
   
   If $r=1$ then $u=d$ and $v=3d=u.$
   
   We thus see that if $v>u$ and $v\neq 3u,$ then  $sup_{n\ge 1}\vert cos(na) -cos(nb)\vert \ge 1+cos\left ( {\pi\over 5}\right )>1.8 > {3\over \sqrt 3},$ which proves (i).
   
   (ii) Assume that $u=v,$ and that  $b \notin \pm a +2\pi\Z.$ There exists $\alpha, \beta \in \{1, \dots, u-1\},$ with $\alpha \neq \beta,$ $\alpha \neq u-\beta$ such that $a \in \pm{2\alpha \pi\over u}+2\pi \Z$ and $b \in \pm{2\beta \pi\over u}+2\pi \Z,$ and $gcd(\alpha,u)=gcd(\beta, u)=1.$ It follows from Bezout's identity that there exists $\gamma \in \Z$ such that $\alpha \gamma -1 \in u\Z.$ If $\beta \gamma\pm1\in u\Z$ then we would have $\alpha \beta \gamma \pm \alpha \in \alpha u\Z\subset u\Z,$ and $\beta \pm\alpha \in u\Z,$ which is impossible. Hence $\gamma \beta -w \in u\Z$ for some $w \in  \{2,\dots, u-2\},$ $gcd(w,u)=1$ since $gcd(\gamma,u)=gcd(\beta, u)=1,$ and we have

$$sup_{n\ge 1}\vert cos(na) -cos(nb)\vert \ge sup_{n\ge 1}\vert \cos(n\gamma a)-cos(n\gamma b)\vert $$ $$=sup_{n\ge 1}\left \vert cos \left ({2n\pi \over u}\right ) -cos \left ( {2nw\pi\over u}\right )\right \vert \ge sup_{n\ge 1}\left \vert cos \left ({2n\alpha \pi \over u}\right ) -cos \left ( {2n\alpha w\pi\over u}\right )\right \vert$$ $$=sup_{n\ge 1}\left \vert cos \left ({2n\alpha \pi \over u}\right ) -cos \left ( {2n\beta \pi\over u}\right )\right \vert = sup_{n\ge 1}\vert cos(na) -cos(nb)\vert.$$

By replacing $w$ by $u-w$ if necessary, we can assume that $2 \le w \le {u\over 2}.$

Now let $w \in \Z$ such that $gcd(u,w)=1.$ We have $a={2\alpha \pi\over u},$ with $gcd(\alpha,u)=1.$ The same argument as above shows that we have

$$  sup_{n\ge 1} \left \vert cos\left ({2 n\pi\over u}\right) -cos \left ( {2nw\pi\over u}\right ) \right \vert =sup_{n\ge 1}\cos(na)-cos(nb)\vert ,$$

with $b={2w\alpha\pi\over u},$ which has order $u.$
   
   (iii) Now assume that $v=3u.$ There exists $\alpha  \in \{1, \dots, u-1\}$ and $\beta \in \{1, \dots, 3u-1\}$ such that $a \in \pm{2\alpha \pi\over u}+2\pi \Z$ and $b \in \pm{2\beta \pi\over 3u}+2\pi \Z,$ and $gcd(\alpha,u)=gcd(\beta, 3u)=1.$ Let $\gamma \in \Z$ such that $\beta \gamma -1 \in 3u\Z.$ Then $gcd(\gamma, 3u)=1,$ and a fortiori $gcd(\gamma,u)=1.$ There exists $w \in \Z$ such that $\alpha \gamma \in \pm w +u\Z,$ and we see as above that we have
   
   $$sup_{n\ge 1} \left \vert cos(na) -cos(nb)\right \vert =sup_{n\ge 1}\left \vert cos \left ({2n\alpha \pi \over u}\right ) -cos \left ( {2n \beta \pi\over 3u}\right )\right \vert$$
   $$=sup_{n\ge 1}\left \vert cos \left ({2n\alpha \gamma \pi \over u}\right ) -cos \left ( {2n \beta \gamma \pi\over 3u}\right )\right \vert=sup_{n\ge 1}\left \vert cos \left ({2nw \pi \over u}\right ) -cos \left ( {2n \pi\over 3u}\right )\right \vert.$$
   
   Conversely let $a ={2 \alpha \pi\over u} \in \pi \Q$ have order $u,$ and let $w\in \Z$ be such that $gcd(u,w)=1.$ If $\alpha$ is not divisible by 3, then $gcd(\alpha,3u)=1.$ If $\alpha$  is divisible by $3,$ then $u$ is not divisible by $3,$ and so $\alpha +u \in \alpha +u\Z$ is not divisible by $3.$ So we can assume without loss of generality that $\alpha$ is not divisible by $3,$ and there exists $\beta \ge 1$ such that $\alpha\beta -1\in 3u\pi\Z.$ Similarly we can assume without loss of generality that $w$ is not divisible by $3,$ and there exists $\gamma \ge 1$ such that $w\gamma -1\in 3u\pi\Z.$  Set $b={2\alpha \gamma \pi\over 3u}.$ Then $b$ has order $3u,$ and we see as above that we have
   
$$sup_{n\ge 1}\left \vert cos \left ({2nw \pi \over u}\right ) -cos \left ( {2n \pi\over 3u}\right )\right \vert \ge   sup_{n\ge 1}\left \vert cos \left ({2n\alpha \gamma w \pi \over u}\right ) -cos \left ( {2n\alpha \gamma \pi\over 3u}\right )\right \vert $$ $$=sup_{n\ge 1}\left \vert cos(na)-cos(nb)\right \vert \ge sup_{n\ge 1}\left \vert cos \left ({2n\alpha \gamma w\beta w \pi \over u}\right ) -cos \left ( {2n\alpha \gamma \beta w\pi\over 3u}\right )\right \vert $$ $$=sup_{n\ge 1}\left \vert cos \left ({2nw \pi \over u}\right ) -cos \left ( {2n \pi\over 3u}\right )\right \vert,$$
   
 which concludes the proof of (iii).

   (iv) Clearly, the first assertion of (iv) is a reformulation of the first assertion of (iii). Now assume that the order $u$ of $a \in \pi \Q$ is divisible by 3, set $v={u\over 3},$ write $a={2\alpha \pi\over u},$ and let $w \in \Z$ such that $gcd(w,v)=1.$  We see as above that we can assume without loss of generality that $gcd(u,w)=1.$

   Since $gcd(\alpha,u)=1,$ we have a fortiori $gcd(\alpha,v)=1,$ so that $gcd(\alpha w, v)=1,$ so that $b:={6\alpha w\over u}$ has order $v$ and we see as above that $a,b, u$ and $w$ satisfy (4). $\square$

    In order to use lemma 3.4, we introduce the following notions.
    
    \begin{definition} Let $u \ge 2,$ and denote by $\Delta(u)$ the set of all integers $s$ satisfying $1 \le s \le {u\over 2},$ $gcd(u,s)=1,$ and let $\Delta_1(u)=\Delta(u)\setminus\{1\}.$ We set
    
        $$\sigma(u)=inf_{w \in \Delta(u)} \left [ sup_{n\ge 1}\left \vert cos\left ({2\pi\over 3u}\right ) -cos \left ({2w\pi\over u}\right )\right \vert\right ],$$
    
    $$\theta(u)=inf_{w \in \Delta_1(u)} \left [ sup_{n\ge 1}\left \vert cos\left ({2\pi\over u}\right ) -cos \left ({2w\pi\over u}\right )\right \vert\right ].$$

    with the convention $\theta(u)=2$ if $\Delta_1(u)=\emptyset.$
    \end{definition}
    
    Notice that $\Delta_1(u)=\emptyset$ if $u=2, 3, 4$ or $6$, and that $\Delta_1(u)\neq \emptyset$ otherwise since as we observed above $H(n) =card \left ( \left (\Z/n\Z\right )^{\times}\right )\ge 3$ if $n\notin \{1,2,3,4,6\}.$
    
  We obtain the following corollary, which shows in particular that the value of $k(a)$ depends only on the order of $a.$
    
    \begin{corollary} Let $a\in \pi\Q,$ and let $u\ge 1$ be the order of $a.$
    
        \smallskip
    
    (i) If $u$ is not divisible by $3,$ then $k(a) =inf(\sigma(u), \theta(u)).$
    
    \smallskip
    
       (ii) If $u$ is divisible by $3,$ then $k(a) =inf(\sigma \left ({u\over 3}\right ), \sigma(u), \theta(u)).$
       
       \end{corollary}
       
       Proof: Set \begin{itemize}
       
       \item $\Lambda_1(a)=\left \{ b\in \pi\Q \ \vert \ b\notin \pm a +2\pi \Z, ord(b)=ord(a)\right \},$
       \item $\Lambda_2(a)= \left \{ b\in \pi\Q \ \vert \  ord(b)=3ord(a) \right \},$ \item $\Lambda_3(a)=\left \{ b\in \pi\Q \ \vert \ 3ord(b)=ord(a) \right \},$ \item $\Lambda_4(a)=\left \{ b\in \pi\Q \ \vert \  ord(b)\neq ord(a)\neq 3ord(b) \right \},$ \end{itemize} and for $1\le i \le 4,$ set
       
       $$\lambda_i(a)=inf_{b\in \Lambda_i(a)}\sup_{n\ge 1}\vert cos(na)-cos(nb)\vert,$$
       
       with the convention $\lambda_i(a)=2$ if $\Lambda_i(a)=\emptyset.$
       
       Since $b\notin \pm a+2\pi\Z$ if $ord(b)\neq ord(a),$ we have $\lambda_2(a)\le {8\over 3\sqrt 3},$ and it follows from lemma 3.4(i) that we have
       
       $$k(a) = inf_{1\le i \le 4}\lambda_i(a)=inf_{1\le i \le 3}\lambda_i(a),$$
       
       and it follows from lemma 3.4 (ii), (iii) and (iv) that $\lambda_1(a)=\theta(u)$ if $\Delta_1(u)\neq \emptyset,$ that $\lambda_{2}(a)=\sigma(u),$ and that  $\lambda_{3}(a)=\sigma\left ({u\over 3} \right )$ if $u$ is divisible by $3.$ $\square$
       
       We have the following result.

   \begin{theorem} Let $m<{8\over 3\sqrt 3}.$ Then the set $\Omega(m):=\{ a \in [0,\pi] : k(a) \le m\}$ is finite.
   \end{theorem}
   
   Proof: It follows from lemma 3.3 applied to ${2\pi\over u}$ and ${6\pi\over u}$ that there exists $u_0\ge 1$ such that we have, for $u\ge u_0,$ 
   
   $$(i) sup_{n\ge 1}\left \vert cos\left ({2n\pi\over u}\right )- cos\left ({2wn\pi\over u}\right )\right \vert >m \ \ \mbox{if} \ 2 \le w \le  inf\left ({u\over 2}, 6\right ),$$
   
    $$(ii) sup_{n\ge 1}\left \vert cos\left ({6n\pi\over u}\right )- cos\left ({2(3w+1)n\pi\over u}\right )\right \vert >m \ \ \mbox{if} \  0 \le w \le 6,$$
    
     $$(iii) sup_{n\ge 1}\left \vert cos\left ({6n\pi\over u}\right )- cos\left ({2(3w+2)n\pi\over u}\right )\right \vert >m \ \ \mbox{if} \  0 \le w \le 6.$$

       Let $u \ge u_0,$ and let $w$ be an integer such that $2 \le w \le {u\over 2}.$ Il ${2w\pi\over u}\le {\pi/2},$ or if ${2w\pi\over u}\ge {5\pi\over 6},$ it follows from  lemma 3.2 and property (i) that we have
       
       $$sup_{n\ge 1}\left \vert cos\left ( {2n\pi\over u}\right ) - cos\left ( {2wn\pi\over u}\right )\right \vert >m.$$
       
       Now assume that ${\pi\over 2} \le {2w\pi\over u}\le{5\pi\over 6}.$ If $\left \vert w -{u\over 3}\right \vert \ge 7,$ it  follows from  lemma 3.2 that we have
       
        $$sup_{n\ge 1}\left \vert cos\left ( {2n\pi\over u}\right ) - cos\left ( {2wn\pi\over u}\right )\right \vert >1.55 >m.$$
        
        If $\left \vert w -{u\over 3}\right \vert < 7,$ set $r =\left \vert 3w-u\right \vert.$ Then $0\le r \le 20,$ and we have
        
        $$sup_{n\ge 1}\left \vert cos\left ( {2n\pi\over u}\right ) - cos\left ( {2wn\pi\over u}\right )\right \vert  \ge sup_{n\ge 1}\left \vert cos\left ( {6n\pi\over u}\right ) - cos\left ( {2nr\pi\over u}\right )\right \vert.$$
        
        If $u$ is not divisible by 3, then either $r=3s+1$ or $r=3s+2,$ with $0\le w \le 6,$ and it follows from (ii) and (iii) that we have
        
  $$sup_{n\ge 1}\left \vert cos\left ( {2n\pi\over u}\right ) - cos\left ( {2wn\pi\over u}\right )\right \vert >m.$$
  
  If $u$ is divisible by $3$ then $r$ is also divisible by 3. Set $v={u\over 3}$ and $s={r\over 3}.$ Then $0\le s \le 6,$ and we have
  
          $$sup_{n\ge 1}\left \vert cos\left ( {2n\pi\over u}\right ) - cos\left ( {2wn\pi\over u}\right )\right \vert  \ge sup_{n\ge 1}\left \vert cos\left ( {2n\pi\over v}\right ) - cos\left ( {2ns\pi\over v}\right )\right \vert.$$
          
          If $s \in  \{ 2, 3, 4, 5, 6\}$ it follows from (i) that we have, if $u\ge 3u_0,$
         
            $$sup_{n\ge 1}\left \vert cos\left ( {2n\pi\over v}\right ) - cos\left ( {2sn\pi\over u}\right )\right \vert >m.$$
            
            Now assume that $s=0.$ If $u \ge 15,$ then $v\ge 5,$ and we have

    $$sup_{n\ge 1}\left \vert cos\left ( {2n\pi\over v}\right ) - cos\left ( {2sn\pi\over u}\right )\right \vert = sup_{n\ge 1}\left \vert cos\left ( {2n\pi\over v}\right ) - 1\right \vert \ge 1+ cos\left ({\pi\over 5}\right )>1.8 >m.$$
    
    Now assume that $s=1.$ We have, with $\epsilon=\pm 1,$
    
   $$ sup_{n\ge 1}\left \vert cos\left ( {2n\pi\over u}\right ) - cos\left ( {2wn\pi\over u}\right )\right \vert = sup_{n\ge 1}\left \vert cos\left ( {2n\pi\over 3v}\right ) - cos\left ( {2n\pi\over 3v}+  {2n\epsilon \pi\over 3}\right )\right \vert$$ $$\ge sup_{n\ge 1}\left \vert cos\left ( {2(3n+1)\pi\over 3v}\right ) - cos\left ( {2(3n+1)\pi\over 3v}+  {2\epsilon \pi\over 3}\right )\right \vert =\sqrt 3 \left \vert sin \left ({2n\pi\over v} +{2 \pi \over 3v} +{\epsilon \pi\over 3}\right ) \right \vert .$$
   
   There exists $p\ge 1$ and $q\in \Z$ such that ${\pi\over 2}-{\pi\over v} \le {2p\pi\over v} +{2\pi \over 3v} +{\epsilon \pi\over 3}+2q\pi\le {\pi\over 2}+{\pi\over v},$ and we obtain, for $u\ge 21,$ $w=v\pm 1,$
   
    $$ sup_{n\ge 1}\left \vert cos\left ( {2n\pi\over u}\right ) - cos\left ( {2wn\pi\over u}\right )\right \vert \ge \sqrt 3 cos\left ({\pi\over v}\right ) \ge \sqrt 3 cos\left ({\pi\over 7}\right ) \ge 1.56 >m.$$
    
    We thus see that if $u \ge u_0$ is not divisible par 3, or if $u\ge max(21, 3u_0)$ is divisible by 3, we have, for $2\le w \le {u\over 2},$
    
    $$ sup_{n\ge 1}\left \vert cos\left ( {2n\pi\over u}\right ) - cos\left ( {2wn\pi\over u}\right )\right \vert\ge m.$$
    
    It follows then  from corollary 3.6 that if the order $u$ of $a \in [0,2\pi]$ satisfies $u \ge max(21, 3u_0),$ we have $k(a) >m.$

      $\square$
      
      We now want to identify the real numbers  $a$ for which $k(a)\le 1.5.$ 
      
      
      
      
      
      
      

       If  $a\in \pi\Q$ has order $1, 2$ or $4,$ then $sup_{n\ge 1} \vert cos(an) -cos(3an)\vert =0.$ We also have the following elementary facts.
      
      \begin{lemma} Let $a\in \pi\Q,$ and let $u\notin\{1,2,4\}$ be the order of $a.$

  \begin{enumerate}
                            \item If $u\notin\{ 3, 5, 6, 8, 9, 10, 11, 12, 15, 16, 18, 22, 24, 30\}$  then 
                            
                            $$sup_{n\ge 1} \vert cos(an) -cos(3an)\vert >1.5.$$
                            
                        \item If $u\in \{3, 6, 9, 12, 15, 18, 24,30\},$ then  $$sup_{n\ge 1} \vert cos(an) -cos(3an)\vert =1.5.$$

                         \item   If $u\in \{5,10\},$ then
                            
                            $$sup_{n\ge 1} \vert cos(an) -cos(3an)\vert ={\sqrt{5}\over 2}.$$
                            
                            \item If $u\in \{8,16\},$ then
                            
                            $$sup_{n\ge 1} \vert cos(an) -cos(3an)\vert =\sqrt 2.$$
                            
                          \item  If $u\in \{11,22\},$ then
                            
                              $$sup_{n\ge 1} \vert cos(an) -cos(3an)\vert =-cos\left ({8\pi\over 11}\right )+cos\left ({24\pi\over 11}\right )=cos\left ({2\pi\over 11}\right )+cos\left ({3\pi\over 11}\right )\approx1.4961.$$
        \end{enumerate}                      
                              
 \end{lemma}
 
 Proof:  We have  $\{e^{ian}\}_{n\ge 1}=\{e^{2in\pi\over u}\}_{1\le n \le u},$ and so we have
 
 $$sup_{n\ge 1} \vert cos(an) -cos(3an)\vert =sup_{n\ge 1} \left \vert cos\left ({2n\pi\over u}\right ) -cos\left ({6n\pi\over u}\right )\right \vert $$ $$=sup_{1\le n\le u} \left \vert cos\left ({2n\pi\over u}\right )- cos\left ({6n\pi\over u}\right )\right \vert ,$$ and the value of $sup_{n\ge 1} \vert cos(an) -cos(3an)\vert $ depends only on the order $u$  of $a.$

 The function $x \to cos(x)-cos(3x)$ is increasing on $\left [0, arccos\left (\frac{1}{\sqrt 3}\right )\right ]$ and decreasing on $\left [arccos\left (\frac{1}{\sqrt 3}\right ), -arccos\left (\frac{1}{\sqrt 3}\right )\right ],$ and $0.275\pi <arccos\left (\frac{1}{\sqrt 3}\right )<0.333\pi.$ Since $cos(x)-cos(3x)>1.5$ if $x=0.275\pi$ or if $x=0.333\pi,$  there exists a closed interval
 $I$ of length $0.058\pi$ on which $ cos (x)-cos(3x)>1.5.$ So if $u \ge 35 >{2\over 0.058},$ there exists $n\ge 1$ such that ${2n\pi\over u}\in I,$ and we have
 
  $$sup_{n\ge 1} \vert cos(an) -cos(3an)\vert >1.5 \ \ \ \forall n\ge 35.$$
  
  The other properties follow from computations of $sup_{1\le n\le u}\left \vert cos\left ({2n\pi\over u}\right ) -cos\left ({6n\pi\over u}\right )\right \vert $ for $3\le u \le 34$ which are left to the reader. $\square$
  
 We now wish to obtain similar estimates for $sup_{n\ge 1}\left \vert cos\left ({2\pi\over n}\right )-cos\left ( {2s\pi\over n}\right )\right \vert$ for $s\in \{2,4,5,6\}.$
 Set  $f_s(x)= cos(x)-cos(sx), \theta_s = sup_{x \ge 0}\vert f(s) \vert, \delta_s=sup_{x \ge 0}\vert f"(s) \vert .$ We have $\theta_s=2$ if $s$ is even, and a computer verification shows that $\theta_s>1.8$ for $s=5.$  It follows from the Taylor-Lagrange inequality that if $f_s$ attains it maximum at $\alpha_s,$ then we have, 

$$\left \vert f_s(x) - \theta_s\right \vert \le {\delta_s\over 2}\left \vert (x-\alpha_s)^2\right \vert, \left \vert f_s(x)\right \vert \ge \theta_s -{\delta_s\over 2}\left \vert (x-\alpha_s)^2\right \vert,$$ and so $\vert f_s(x)\vert >1.5$ if $\left \vert (x-\alpha_s)^2\right \vert \le {2\theta_s-3\over \delta_s}.$ So if $l_s <\sqrt{2\theta_s-3\over \delta_s}$ there exists a closed interval of length $2l_s$ on which $\vert f_s(x)\vert >1.5.$ Let $u_s \ge {\pi\over l_s}$ be an integer. We obtain

\begin{equation} sup_{n\ge 1}\left \vert cos \left ( {2n\pi\over u}\right ) -cos \left ( {2sn\pi\over u}\right )\right \vert  >  1.5 \ \ \forall u\ge u_s.\end{equation}

 Values for $u_s$ are given by the following table.

$$\begin{array}{| l | l | l | c | r |}\hline  s & \theta_s & \delta_s & l_s & u_s \cr \hline 2 & 2 &\le 5 &0.4472 & 8\cr \hline
 4 & 2 &\le 17 &0.2425 & 13\cr \hline 5 & >1.8 &\le 26 &0.1519 & 21\cr \hline  6 & 2 &\le 37 &0.1644 &20 \cr \hline   \end{array}$$
 
 We obtain the following result
 
 \begin{lemma} Let $u\ge 4$ be an integer, and let $s\le {u\over 4}$ be a nonnegative integer.
 
 If $s \neq 3,$ then we have
 
 $$sup_{n\ge 1}\left \vert cos \left ( {2n\pi\over u} \right )-  cos \left ( {2ns\pi\over u} \right )\right \vert >1.5$$
 
 \end{lemma}
 
 Proof: If $s=0,$ then 
 
  $$sup_{n\ge 1}\left \vert cos \left ( {2n\pi\over u} \right )-  cos \left ( {2ns\pi\over u} \right )\right \vert =sup_{n\ge 1}\left \vert cos \left ( {2n\pi\over u} \right )- 1 \right \vert  >1.8.$$
  
  If $s\ge 7,$ the result follows from lemma 3.2 (i). If $s\in \{2, 4,6\},$ the result follows from the table since $u \ge 4s.$ If $s=5,$ the result also follows from the table for $u \ge 21,$ and a direct computation shows that we have
  
    $$sup_{n\ge 1}\left \vert cos \left ( {2n\pi\over 20} \right )-  cos \left ( {10n\pi\over 20} \right )\right \vert =sup_{1\le n \le 20}\left \vert cos \left ( {n\pi\over 10} \right )-  cos \left ( {n\pi\over 2} \right )\right \vert =1+cos \left ( {\pi\over 5}\right )>1.80.$$ $\square$

Now set  $g_s(x)= cos(3x)-cos(sx), \theta_s = sup_{x \ge 0}\vert g(s) \vert, \delta_s=sup_{x \ge 0}\vert g"(s) \vert .$ We have $\theta_s=2$ if $s$ is even, and a computer verification shows that $\theta_s>1.85$ for $s=5,$ $\theta_s>1.91$ for $s=7, s=11,$ $\theta_s>1.97$ for $s=13,s=17,$ $\theta _s>1.96$ for $s=19.$ We see as above that if $l_s <\sqrt{2\theta_s-3\over \delta_s},$ and if  $u_s \ge {\pi\over l_s}$ is an integer, we have

\begin{equation} sup_{n\ge 1}\left \vert cos \left ( {2sn\pi\over u}\right ) -cos \left ( {6n\pi\over u}\right )\right  \vert >1.5 \ \ \forall u\ge u_s.\end{equation}

 We have the following table.

$$\begin{array}{| l | l | l | c | r |}\hline  s & \theta_s & \delta_s & l_s & u_s \cr \hline 2 & 2 &\le 13 &0.2774 & 12\cr \hline 4 & 2 &\le 23 &0.2085 & 16\cr \hline 
 5 & > 1.85 &\le 34 &0.1435 & 22\cr \hline 7 & >1.91 &\le 58 &0.1189 & 27\cr \hline 8 & 2 &\le 73 &0.1170 &27 \cr \hline  10 & 2 &\le 109 &0.0958 &33 \cr \hline 11 & > 1.91 &\le 130 &0.0794 &40 \cr \hline 
 13 & > 1.97 & \le 178 &0.0727 & 44\cr \hline 14 & 2 &\le 205 &0.0698 & 45\cr \hline  16 & 2  &\le 275 &0.0603 &53 \cr \hline  17 & > 1.97 &\le 298 &0.0562 &56 \cr \hline 19 & > 1.96 &\le 390 &0.0486 & 65\cr \hline 
 20 & 2 & \le 409 &0.0494 &64 \cr \hline  \end{array}$$
We will be interested here to the case where $u$ is not divisible by $3$ and where ${2s\pi\over u}\le {\pi\over 2},$ which means that $ u \ge 4s.$ So we are left with $s=2,$ $u=8,10$ or $11,$ and with $s=5,$ $u=20.$  We obtain, by direct computations

$$sup_{n\ge 1}\left \vert cos \left ( {4n\pi\over 8}\right ) -cos\left ( {6n\pi\over 8}\right )\right\vert =sup_{n\ge 1}\left \vert cos \left ( {n\pi\over 2}\right ) -cos\left ( {3n\pi\over 4}\right )\right \vert =2.$$


$$sup_{n\ge 1}\left \vert cos \left ( {4n\pi\over 10}\right ) -cos\left ( {6n\pi\over 10}\right )\right\vert =sup_{n\ge 1}\left \vert cos \left ( {2n\pi\over 5}\right ) -cos\left ( {3n\pi\over 5}\right )\right\vert =2.$$

$$sup_{n\ge 1}\left \vert cos \left ( {4n\pi\over 11}\right ) -cos\left ( {6n\pi\over 11}\right )\right\vert  =cos\left ({20\pi\over 11}\right ) -cos\left ({30\pi\over 11}\right ) = cos\left ({2\pi\over 11}\right ) +cos\left ({3\pi\over 11}\right )\approx 1.4961.$$

$$sup_{n\ge 1}\left \vert cos \left ( {10n\pi\over 20}\right ) -cos\left ( {6n\pi\over 20}\right )\right\vert =sup_{n\ge 1}\left \vert cos \left ( {n\pi\over 2}\right ) -cos\left ( {3n\pi\over 10}\right )\right \vert >1.80.$$


We obtain the following lemma.

\begin{lemma} Let $u,s$ be positive integers satisfying $u\ge 4,$ ${u\over 4}\le s\le {5u\over 12},$ with $s\ge 2,$ so that $u\ge 5.$
We have
$$sup_{n\ge 1}\left \vert cos \left ( {2n\pi\over u}\right ) -cos\left ( {2sn\pi\over u}\right )\right\vert  \left \{ \begin{array}{l} =cos\left ({\pi\over 5}\right )+ cos\left ({2\pi\over 5}\right )\ \mbox{if} \ u=5, s=2, \ \mbox{or if}\ u=10, s=3,\cr= \sqrt 2 \ \mbox{if} \ u=8, s=3, \ \mbox{or if}\ u=16, s=5, \cr =cos\left ({2\pi\over 11}\right )+ cos\left ({3\pi\over 11}\right ) \   \mbox{if} \ u=11, s=4 \ \mbox{or if} \ u=22, s=7,\\ =1.5 \ \mbox{if} \ u=12, s=3,\\ >1.5 \ \mbox{otherwise}.\end{array} \right .$$



\end{lemma}

Proof: Set $r=\vert 3s-u\vert.$ Since ${2\pi\over 3} -{\pi\over 2}={5\pi\over 6}-{2\pi\over 3}={\pi\over 6},$ we have $0\le{2\pi r \over u}\le{\pi\over 2}.$ If $r\ge 21,$ it follows from lemma 3.1(i) that $sup_{n\ge 1}\left \vert cos \left ( {2n\pi\over u}\right ) -cos\left ( {2sn\pi\over u}\right )\right\vert >1.5.$

If $u$ is not divisible by $3,$ then $r$ is not divisible by $3$ either, and it follows from the discussion above that if $r\neq 1$ and $r\neq 2,$ we have $$sup_{n\ge 1}\left \vert cos \left ( {2n\pi\over u}\right ) -cos\left ( {2sn\pi\over u}\right )\right\vert \ge sup_{n\ge 1}\left \vert cos\left ({6n\pi\over u}\right )- cos\left ( {2rn\pi\over u}\right )\right \vert >1.5.$$ 

The condition $r=2$ gives $\left \vert s -{u\over 3}\right \vert ={2\over 3}.$ We saw above that in this situation $sup_{n\ge 1}\left \vert cos \left ( {2n\pi\over u}\right ) -cos\left ( {2sn\pi\over u}\right )\right\vert >1.5$ unless $u=11,$ which gives $s=3,$ and we have

$$sup_{n\ge 1}\left \vert cos\left ({2n\pi\over 11}\right ) -cos\left ({6n\pi\over 11}\right ) \right \vert =sup_{1\le n \le 11}\left \vert cos\left ({6n\pi\over 11}\right ) -cos\left ({4n\pi\over 11}\right ) \right \vert $$ $$= \left \vert cos\left ({30\pi\over 11}\right ) +cos\left ({20\pi\over 11}\right )\right \vert =cos\left ({2\pi\over 11}\right ) +cos\left ({3\pi\over 11}\right )\approx 1.4961.$$

The condition $r=1$ gives $\left \vert s -{u\over 3}\right \vert ={1\over 3},$ which gives $s={u-1\over 3}$ if $u\equiv 1$ mod 3, and $s={u+1\over 3}$ if $u\equiv 2$ mod 3. In this situation we have

$$sup_{n\ge 1}\left \vert cos \left ( {2n\pi\over u}\right ) -cos\left ( {2sn\pi\over u}\right )\right\vert \ge sup_{n\ge 1}\left \vert cos \left ( {2n\pi\over u}\right ) -cos\left ( {6n\pi\over u}\right )\right\vert.$$

Since we must have $\left \vert s -{u\over 3}\right \vert ={1\over 3},$ it follows from lemma 3.8 that if $n\notin \{5, 8,10,11, 16, 22\}, $ or if $u=5, s\neq 2,$ or if $u=8,s\neq 3,$ or if $u=10, s\neq 3,$ or if $u=11,$ $s\neq 4,$ or if $u=16,$ $s\neq 5,$ or if $u=22, s\neq 7$ we have

$$sup_{n\ge 1}\left \vert cos \left ( {2n\pi\over u}\right ) -cos\left ( {2sn\pi\over u}\right )\right\vert >1.5.$$

A direct computation shows the that we have

$$sup_{n\ge 1}\left \vert cos \left ( {2n\pi\over u}\right ) -cos\left ( {2sn\pi\over u}\right )\right\vert =sup_{1\le n\le u}\left \vert cos \left ( {2n\pi\over u}\right ) -cos\left ( {2sn\pi\over u}\right )\right\vert$$ $$= \left \{ \begin{array}{l} cos\left ({\pi\over 5}\right )+ cos\left ({2\pi\over 5}\right )\ \mbox{if} \ u=5, s=2, \ \mbox{or if}\ u=10, s=3,\cr \sqrt 2 \ \mbox{if} \ u=8, s=3, \ \mbox{or if}\ u=16, s=5, \cr cos\left ({2\pi\over 11}\right )+ cos\left ({3\pi\over 11}\right ) \ \mbox{if} \ u=11, s=4 \ \mbox{ or if} \ u=22, s=7.\end{array} \right .$$

We now consider the case where $u=3v$ is divisible by $3.$ Then $r$ is also divisible by 3. If $r=0,$ and if $u\neq 9,$ then we have 

$$sup_{n\ge 1}\left \vert cos \left ( {2n\pi\over u}\right ) -cos\left ( {2sn\pi\over u}\right )\right\vert \ge sup_{n\ge 1}\left \vert cos \left ( {2n\pi\over v}\right ) -1\right\vert>1.8.$$

If $u=9,$ then $s=3,$ and we have

$$sup_{n\ge 1}\left \vert cos \left ( {2n\pi\over u}\right ) -cos\left ( {2sn\pi\over u}\right )\right\vert =sup_{1\le n\le 9}\left \vert cos \left ( {2n\pi\over 9}\right ) -cos\left ( {2n\pi\over 3}\right )\right\vert =1.5.$$

Now assume that $r=3,$ which means that $s=v+\epsilon,$ with $\epsilon=\pm1.$ We have

$$sup_{n\ge 1}\left \vert cos \left ( {2n\pi\over u}\right ) -cos\left ( {2sn\pi\over u}\right )\right\vert =sup_{1\le n \le 3v}\left \vert cos \left ( {2n\pi\over 3v}\right ) -cos\left ( {2n\pi\over 3}+{2\epsilon n\pi\over 3v}\right )\right\vert$$ $$=2sup_{1\le n\le 3v}\left \vert \sin \left ( {n\pi\over 3}+{(1+\epsilon)n\pi\over 3v} \right ) \right \vert \left \vert sin \left ( -{n\pi\over 3}+{(1-\epsilon)n\pi\over 3v} \right ) \right \vert$$ $$ =sup_{1\le n \le 3v}\left \vert sin\left ({n\pi\over 3}\right )\right \vert \left \vert  sin\left ({n\pi\over 3} +{2n\pi\over 3v}\right ) \right \vert $$ $$\ge\sqrt 3 sup_{0\le n \le v}\left \vert sin\left ({(3n+1)\pi\over 3}+{2(3n+1)\pi\over 3v} \right ) \right \vert$$ $$=\sqrt 3 sup_{0\le n \le v}\left \vert sin\left ({2n\pi\over v} + {(
v+2)\pi\over 3v}\right ) \right \vert.$$

Since $sin(x) >{\sqrt 3\over 2}$ for ${\pi\over 3}<x<{2\pi\over 3},$ there exists $n\in \{1,\dots,v\}$ such that $ sin\left ({2n\pi\over v} + {(v+2)\pi\over 3v}\right ) >{\sqrt 3\over 2}$ if $v\ge 7,$ and we obtain

$$sup_{n\ge 1}\left \vert cos \left ( {2n\pi\over u}\right ) -cos\left ( {2sn\pi\over u}\right )\right\vert  >1.5\ \ \mbox{if} \ u\ge 21.$$

 We are left with the cases where $u=6, v=2, s=1$ or $3$, $u=9, v=3, s=2$ or $4$, $u=12, v=4, s=3$ or $5$, $u=15, v=5, s=4$ or $6$, $u=18, v=6, s=5$ or $7.$ But $s=1$ is not relevant, and the condition ${u\over 4} \le s \le {5u\over 12}$ is not satisfied for $u=6, s=3$ and for $u=9, s=2$ or $4.$

      Direct computations which are left to the reader show that we have

   $$ sup_{n\ge 1}\left \vert cos \left ( {2n\pi\over u}\right ) -cos\left ( {2sn\pi\over u}\right )\right\vert \left \{\begin{array}{l}  >1.64 
\  \mbox{if} \ u=15 \ \mbox{and} \ s=4,\\ > 1.70  \ \mbox{or if}
    \ u=18\ \mbox{and} \ s=5\ \mbox{or} \ s=7\\  > 1.72  \ \mbox{ if}
    \ u=15\ \mbox{and} \ s=6,\\  > 1.73  \ \mbox{ if}
    \ u=12\ \mbox{and} \ s=5.\\\end{array}\right .$$

 So $sup_{n\ge 1}\left \vert cos \left ( {2n\pi\over u}\right ) -cos\left ( {2sn\pi\over u}\right )\right\vert  >1.5$ if ${u\over 4} \le s \le {5u\over 12}$ when $u$ is divisible by $3$ and when $s-{u \over 3}\in \{-1,0,1\},$ unless $u=12$ and $s=3.$ If $u=12$ and $s=3,$ we have $$ sup_{n\ge 1}\left \vert cos \left ( {2n\pi\over u}\right ) -cos\left ( {2sn\pi\over u}\right )\right\vert= sup_{n\ge 1}\left \vert cos \left ( {n\pi\over 6}\right ) -cos\left ( {n\pi\over 2}\right )\right\vert=1.5.$$

 Now assume that $u=3v$ is divisible by $3,$ and that $2\le \left \vert s - v \right \vert\le 6.$ Set again $r=\left \vert 3s-u\right \vert,$ and set $p={r\over 3},$ so that $2\le p \le 6.$ Notice also that $p\le {u\over 12}$ since $r\le {u\over 4},$ so that $u\ge 24$ and $v\ge 8.$ We have
 
 $$sup_{n\ge 1}\left \vert cos \left ( {2n\pi\over u}\right ) -cos\left ( {2sn\pi\over u}\right )\right\vert\ge sup_{n\ge 1}\left \vert cos \left ( {6n\pi\over u}\right ) -cos\left ( {2rn\pi\over u}\right )\right\vert$$ $$=sup_{n\ge 1}\left \vert cos \left ( {2n\pi\over v}\right ) -cos\left ( {2pn\pi\over u}\right )\right\vert.$$
 
 It follows then from lemma 3.9 that $sup_{n\ge 1}\left \vert cos \left ( {2n\pi\over u}\right ) -cos\left ( {2sn\pi\over u}\right )\right\vert>1.5$ if $p\neq 3.$
 
 If $p=3,$ then $u\ge 36,$ and so $v\ge 12.$ Since $s-v=\pm 3,$ it follows from lemma 3.5 that we only have to consider the following cases:
 
 \begin{itemize}
 
 \item $u=36, s =9$ or $15,$
 
  \item $u=45, s =12$ or $18,$
  
   \item $u=54, s =15$ or $21,$
   
    \item $u=72, s =21$ or $27,$
    
     \item $u=90, s =27$ or $33.$
     
     Direct computations which are left to the reader show that we have

     \end{itemize}

   $$ sup_{n\ge 1}\left \vert cos \left ( {2n\pi\over u}\right ) -cos\left ( {2sn\pi\over u}\right )\right\vert \left \{\begin{array}{l} >1.93 \ \mbox{if} \ u=36 \ \mbox{and} \ s=9,  \ \mbox{or if}
 \  u=45\ \mbox{and} \ s=12 \ \mbox{or} \ 18,\\  \ \mbox{or if} \ u=72 \ \mbox{and} \ s=27, \ \mbox{or if} \ u=90 \ \mbox{and} \ s=27 \ \mbox{or} \ 33,\\ > 1.91  \ \mbox{or if}
    \ u=54\ \mbox{and} \ s=15,\\  > 1.87  \ \mbox{or if}
    \ u=72\ \mbox{and} \ s=24,\\  > 1.85  \ \mbox{or if}
    \ u=36\ \mbox{and} \ s=15,\\  > 1.83  \ \mbox{or if}
    \ u=54\ \mbox{and} \ s=21.\\\end{array}\right .$$
      
      This concludes the proof of the lemma. $\square$

      \begin{lemma} Let $u,s$ be positive integers satisfying $ {5u\over 12} \le s \le {u\over 2},$ with $s\ge 2,$ so that $u\ge 4.$
We have
 $$ sup_{n\ge 1}\left \vert cos \left ( {2n\pi\over u}\right ) -cos\left ( {2sn\pi\over u}\right )\right\vert=\left \{\begin{array}{l} =1.5 \ \mbox{if} \ u=6\ \mbox{and} \ s=3, \\ >1.5\ \mbox{ otherwise}
  \end{array}\right .$$

\end{lemma}

Proof: If $s\ge 4,$ it follows from lemma 3.2(ii) that we have

$$ sup_{n\ge 1}\left \vert cos \left ( {2n\pi\over u}\right ) -cos\left ( {2sn\pi\over u}\right )\right\vert  >1.57.$$

So we only have to consider the cases $s=3,$ $u=6$ or $7,$ and $s=2,$ $u=4.$

A direct computation then shows that we have

   $$ sup_{n\ge 1}\left \vert cos \left ( {2n\pi\over u}\right ) -cos\left ( {2sn\pi\over u}\right )\right\vert \left \{\begin{array}{l} =2 \ \mbox{if} \ u=4 \ \mbox{and} \ s=2, \\ =1.5  \ \mbox{if}
    \ u=6\ \mbox{and} \ s=3,\\  =cos\left ({2\pi\over 7}\right ) + cos\left ({\pi\over 7}\right )\approx 1.5245 \ \mbox{ if}
    \ u=7\ \mbox{and} \ s=3.\\\end{array}\right .$$
    
    $\square$
    
    We consider again the numbers $\theta(u)$ and $\sigma(u)$ introduced in definition 3.6.

   It follows from lemma 3.8, lemma 3.9, lemma 3.10 and lemma 3.11 that we have the following results.
   
   \begin{lemma} We have \ $\theta (5)=\theta(10)= cos\left ({\pi\over 5}\right )+cos\left ({2\pi\over 5}\right ), \theta(8)=\theta(16)=\sqrt 2, \theta (11)=\theta(22)=cos\left ({2\pi\over 11}\right )+cos\left ({3\pi\over 11}\right ),$ and $\theta(u)>1.5$ for $u\ge 4, u \neq 5, u\neq 8, u \neq 10, u \neq 11, u \neq 16, u \neq 22.$
   
   \end{lemma}
   
      \begin{lemma} We have $\sigma(u)=1.5$ if $u\in \left \{1,2 ,3 ,4 ,5 ,6 ,8,10 \right \},$ $\sigma(u)>1.5$ otherwise.
   
   \end{lemma}
   
   Hence if $u$ is divisible by $3,$ we have $\sigma \left ({u\over 3} \right )=1.5$ if $u\in \left \{3, 6, 9, 12 , 15 , 18, 24,30 \right \},$ $\sigma(u)>1.5$ otherwise. We then deduce from corollary 3.6 a complete description of the set $\Omega(1.5)=\{ a \in [0, \pi] \ \vert \ k(a)\le 1.5\}.$
   
   \begin{theorem} Let $a \in [0,\pi].$

    \begin{itemize}
   
   \item If $a \in \left \{ {\pi\over 5}, {2\pi\over 5}, {3\pi\over 5}, {4\pi\over 5} \right \},$  then $k(a)= cos\left ({\pi\over 5}\right)+  cos\left ({2\pi\over 5}\right)\approx 1,1180.$
   
   \smallskip
   
      \item If $a \in \left \{ {\pi\over 8}, {\pi\over 4}, {3\pi\over 8}, {5\pi\over 8}, {5\pi\over 4}, {7\pi\over 8}\right \},$  then $k(a)= \sqrt 2 \approx 1,4142.$
      
        \smallskip
   
      \item If $a \in \left \{ {\pi\over 11}, {2\pi\over 11}, {3\pi\over 11}, {4\pi\over 11}, {5\pi\over 11}, {6\pi\over 11}, {7\pi\over 11}, {8\pi\over 11}, {9\pi\over 11}, {10\pi\over 11}\right \},$  then $k(a)= cos\left ({2\pi\over 11}\right)+  cos\left ({3\pi\over 11}\right)\approx 1,4961.$
      
           \smallskip
      
      \item If $a \in \left \{ {0}, {\pi\over 6}, {\pi\over 3}, {\pi\over 2}, {2\pi\over 3}, {5\pi\over 6}\right \} \cup  \left \{ {\pi\over 9}, {2\pi\over 9}, {4\pi\over 9}, {5\pi\over 9}, {7\pi\over 9}, {8\pi\over 9}\right \} \cup   \left \{ {\pi\over 12}, {5\pi\over 12}, {7\pi\over 12}\right \}\\ \cup\left \{ {\pi\over 15}, {2\pi\over 15}, {4\pi\over 15}, {7\pi\over 15}, {8\pi\over 15}, {11\pi\over 15}, {13\pi\over 15}, {14\pi\over 15}\right \},$ then $k(a)=1.5.$
      
           \smallskip
           
           \item For all other values of $a,$ we have $1.5 < k(a) \le{8\over 3\sqrt 3}\approx 1.5396.$
      
      \end{itemize}
      
      \end{theorem}
      
      \begin{corollary} Let $G$ be an abelian group, and let $(C(g))_{g\in G}$ be a $G$-cosine family in a unital Banach algebra $A$ such that $sup_{g\in G}\Vert C(g)-c(g)\Vert < 
      {\sqrt 5\over 2}$ for some bounded scalar $G$-cosine family $(c(g))_{g\in G}.$ Then $C(g)=c(g)$ for every $g\in G.$
      \end{corollary}
      
      Proof: Let $g \in G.$ Since the scalar cosine sequence $(c(ng))_{n\in \Z}$ is bounded, a standard argument shows that there  exists  $a(g) \in \R$ such that $c(ng)=cos(na(g))1_A$ for $n\in \Z.$ Since $k(a(g))\ge {\sqrt 5\over 2},$ it follows from corollary 2.4 that $C(ng)=cos(na(g))1_A=c(ng)$ for $n\in \Z,$ and $C(g)=c(g).$ $\square$

{\it IMB, UMR 5251

Universit\'e de Bordeaux

351, cours de la Lib\'eration

33405 Talence (France)

\end{document}